\documentclass[a4,leqno,12pt]{article} 

\usepackage{amsmath,amsthm,amssymb,bm} 
\numberwithin{equation}{section}

\title{{\bf 
Asymptotics of abelian group-partitions \\ 
and associated Dirichlet series}} 
\author{Tetsuya MOMOTANI} 

\newtheorem{The}{Theorem}[section]
\newtheorem{Pro}[The]{Proposition}
\newtheorem{Lem}[The]{Lemma}
\newtheorem{Def}[The]{Definition}

\newtheorem{Cor}[The]{Corollary}

\begin{document} 

%%%%%%%%%%%%%%%%%%%%%%%%%%%%%%%%%%%%%%%%%%%%%%%%%%%%%%%%%%%%%%%%%%%%%%%%%%%%%%%%
\maketitle 

\begin{abstract} 
We introduce a notion of a group-partition for a finite Abelian group, 
which is a generalized notion of the standard partition. 
To obtain asymptotic distributions of group-partition, 
we study the Dirichlet series for group-partitions 
by employing the generating function of the plane partition. 

\end{abstract} 

%%%%%%%%%%%%%%%%%%%%%%%%%%%%%%%%%%%%%%%%%%%%%%%%%%%%%%%%%%%%%%%%%%%%%%%%%%%%%%%%
%%%%%%%%%%%%%%%%%%%%%%%%%%%%%%%%%%%%%%%%%%%%%%%%%%%%%%%%%%%%%%%%%%%%%%%%%%%%%%%%

%%%%%%%%%%%%%%%%%%%%%%%%%%%%%%%%%%%%%%%%%%%%%%%%%%%%%%%%%%%%%%%%%%%%%%%%%%%%%%%%

\section{Introduction} 

Let $ a(n) $ 
be the number of isomorphic classes of Abelian groups of order $n$ 
It is known (Cf. \cite{Sa}) that 
\begin{align} 
\sum_{ 
\substack{ n \leq x  \\ 
n \equiv k \;( \bmod \, j ) } } 
a (n) 
= C(j) x  + o(x) \qquad \text{as } x \to \infty .
\label{AAB}
\end{align} 
Here $j$ and $k$ are positive integers such that $ (j,k) = 1 $
and $ C(j) $ is a positive constant. 

The aim of this paper is a generalization of (\ref{AAB}) 
by introducing a new notion which we call a group partition; 
for an isomorphic class $G$ of finite Abelian groups,  
a sequence $ ( G_{j} ) $ of isomorphic classes of finite Abelian groups 
is called a group-partition of $G$
if $ ( G_{j} ) $ satisfies 
\begin{align*} 
G = G_{1} \oplus G_{2} \oplus G_{3} \oplus \cdots \oplus G_{r}, \quad 
G_{1} \supseteq G_{2} \supseteq G_{3} \supseteq \cdots \supseteq G_{r} .
\end{align*} 
Let $ \pi_{r} (G) $ be the number of group-partitions of $G$ 
with at most $r$ factors. 
Let $ a_{r} (n) $ be an arithmetic function defined by 
$ a_{r} (n) := \sum_{ |G| = n } \pi_{r} (G) $ 
where $G$ runs over isomorphic classes of finite Abelian groups of order $n$. 
Note that $ a_{1} (n) = a (n) $. 
Then the main result claims that 
\begin{align} 
\sum_{ 
\substack{ n \leq x  \\ 
n \equiv k \;( \bmod \, j ) } } 
a_{r} (n)
= C_{r} (j) x  + o(x) \qquad \text{as } x \to \infty . 
\label{ABC}
\end{align} 
Here $ C_{r} (j) $ is a positive constant 
explicitly expressed in terms of 
the Riemann zeta function and the Euler function (Theorem 4.3). 

To prove (\ref{ABC}), 
we introduce the Dirichlet series defined by 
\begin{align*} 
J_{r} (s , \chi ) 
:= \sum_{ n = 1}^{ \infty } 
\chi (n) a_{r} (n) n^{-s} \qquad \mathrm{Re} (s) > 1 ,
\end{align*}
where $ \chi $ is a Dirichlet character. 
Using plane partitions and its generating functions, 
we show that $a_{r} (n) $ is multiplicative in the sence of 
$ a_{r} (mn) = a_{r} (m) a_{r} (n) $ 
when $ (m , n) = 1 $. 
From this function, 
we see that $ J_{r} (s, \chi ) $ has 
an Euler product expression and an analytic continuation to 
$ \mathrm{Re} (s) > 0 $. 
Then the asymptotic distribution (\ref{ABC}) follows 
from the Ikehara-Wiener theorem by the standard technique. 

%%%%%%%%%%%%%%%%%%%%%%%%%%%%%%%%%%%%%%%%%%%%%%%%%%%%%%%%%%%%%%%%%%%%%%%%%%%%%%%%
%%%%%%%%%%%%%%%%%%%%%%%%%%%%%%%%%%%%%%%%%%%%%%%%%%%%%%%%%%%%%%%%%%%%%%%%%%%%%%%%

\section{Group-partitions and plane partitions}

Let $n$ be a positive integer 
and $ ( n_{j} ) $ be a sequence of non-negative integers such that
\begin{align*} 
n = n_{1} + n_{2} + n_{3} + \cdots , \quad 
n_{1} \geq n_{2} \geq n_{3} \geq \cdots .
\end{align*} 
Then $ (n_{j} ) $ is called a partition of $n$, 
and we write $ (n_{j}) \vdash n$. 
We denote by $ P(n) $ the number of all partitions of $n$. 

As a generalization of notion of partitions for positive integers, 
we define group-partitions for finite Abelian groups. 

\begin{Def} 
Let $G$ be an isomorphic class of finite Abelian groups. 
Suppose that $ ( G_{j} ) $ is 
a sequence of isomorphic classes of finite Abelian groups such that 
\begin{align*} 
G = G_{1} \oplus G_{2} \oplus G_{3} \oplus \cdots , \quad 
G_{1} \supseteq G_{2} \supseteq G_{3} \supseteq \cdots .
\end{align*} 
Then we call $ ( G_{j} ) $ a group-partition of $G$, 
and write $ ( G_{j} ) \vdash G $. 
For a positive integer $r$, 
we define $ \pi_{r} (G) $ as the number of group-partitions of $G$ 
with at most $r$ factors. 
We also define $ \pi _{\infty} (G) $ 
as the number of all group-partitions of $G$. \qed
\end{Def} 

We give some examples. \\

\noindent 
{\bf Example 1.}
Suppose that 
$ G = ( \mathbb{Z} / p \mathbb{Z} )^{n} $ 
with a prime number $p$ and a positive integer $n$. 
We put $ G_{j} = ( \mathbb{Z} / p \mathbb{Z} )^{ n_{j} } $. 
Then $ (G_{j}) \vdash G \Leftrightarrow (n_{j}) \vdash n $. 
Hence we have $ \pi_{\infty} (G) = P(n)$. 
The decomposition theorem of finite Abelian groups states that 
every finite Abelian group can be expressed as the direct sum of 
$  ( \mathbb{Z} / p^{m} \mathbb{Z} )^{n} $'s. 
Hence the group-partition is 
a generalization of the usual partition.  \\

\noindent
{\bf Example 2.}
If $ G = (\mathbb{Z} / p \mathbb{Z} )^{2} 
\oplus (\mathbb{Z} / p^{2} \mathbb{Z} ) $ with a prime number $p$, 
then $ \pi_{1} (G) = 1$, $ \pi _{2} (G) = 2 $, 
$ \pi _{3} (G) = \pi _{4} (G) = \cdots = \pi _{ \infty } (G) = 3 $. 
In fact, all of group-partitions of $G$ are given as follows.  

(1) $ G = G_{1} $, where
$ G_{1} = (\mathbb{Z} / p \mathbb{Z} )^{2} 
\oplus (\mathbb{Z} / p^{2} \mathbb{Z} ) $. 

(2) $ G = G_{1} \oplus G_{2} $, 
$ \quad G_{1} = (\mathbb{Z} / p \mathbb{Z} ) 
\oplus (\mathbb{Z} / p^{2} \mathbb{Z} ) $, 
$ G_{2} = (\mathbb{Z} / p \mathbb{Z} )  $. 

(3) $ G = G_{1} \oplus G_{2} \oplus G_{3} $, 
$ \quad G_{1} = (\mathbb{Z} / p^{2} \mathbb{Z} ) $, 
$ G_{2} = (\mathbb{Z} / p \mathbb{Z} ) $, 
$ G_{3} = (\mathbb{Z} / p \mathbb{Z} ) $. \\

\noindent
{\bf Example 3.}
If $ G = (\mathbb{Z} / p \mathbb{Z} )^{3} 
\oplus (\mathbb{Z} / q \mathbb{Z} )^{2} $ 
with prime numbers $p \neq q $, 
then $ \pi_{1} (G) = 1 $, $ \pi_{2} (G) = 4 $, 
$ \pi _{3} (G) = \pi _{4} (G) = \cdots = \pi _{ \infty } (G) = 6 $. 

(1) $ G = G_{1} $, where
$ G_{1} = (\mathbb{Z} / p \mathbb{Z} )^{3} 
\oplus (\mathbb{Z} / q \mathbb{Z} )^{2} $. 

(2) $ G = G_{1} \oplus G_{2} $, 
$ \quad G_{1} = (\mathbb{Z} / p \mathbb{Z} )^{3} 
\oplus (\mathbb{Z} / q \mathbb{Z} ) $, 
$ G_{2} = (\mathbb{Z} / q \mathbb{Z} )  $. 

(3) $ G = G_{1} \oplus G_{2} $, 
$ \quad G_{1} = 
( \mathbb{Z} / p \mathbb{Z} )^{2} 
\oplus (\mathbb{Z} / q \mathbb{Z} )^{2} $, 
$ G_{2} = (\mathbb{Z} / p \mathbb{Z} )  $, 

(4) $ G = G_{1} \oplus G_{2} $, 
$ \quad G_{1} = 
( \mathbb{Z} / p \mathbb{Z} )^{2} 
\oplus (\mathbb{Z} / q \mathbb{Z} ) $, 
$ G_{2} = (\mathbb{Z} / p \mathbb{Z} ) 
\oplus (\mathbb{Z} / q \mathbb{Z} )  $, 

(5) $ G = G_{1} \oplus G_{2} \oplus G_{3} $, 
$ G_{1} = (\mathbb{Z} / p \mathbb{Z} ) 
\oplus (\mathbb{Z} / q \mathbb{Z} ) ^{2} $, 
$ G_{2} = (\mathbb{Z} / p \mathbb{Z} ) $, 
$ G_{3} = (\mathbb{Z} / p \mathbb{Z} ) $. 

(6) $ G = G_{1} \oplus G_{2} \oplus G_{3} $, 
$ G_{1} = (\mathbb{Z} / p \mathbb{Z} ) 
\oplus (\mathbb{Z} / q \mathbb{Z} )  $, 
$ G_{2} = (\mathbb{Z} / p \mathbb{Z} ) 
\oplus (\mathbb{Z} / q \mathbb{Z} )  $, 
$ G_{3} = (\mathbb{Z} / p \mathbb{Z} ) $. \\

In general, it is easy to see the following. 

\begin{Lem} 
Suppose that $ p_{1} , p_{2} , p_{3} , \cdots  $ 
are distinct prime numbers. 
If $ G_{p_{j}} $ is 
an Abelian $ p_{j}$-group, 
then we have
$ \pi_{r} ( G_{ p_{1} } \oplus G_{ p_{2} } \oplus G_{ p_{3}} \oplus  \cdots ) 
= \pi_{r} ( G_{p_{1}} ) \times \pi_{r} (G_{ p_{2} }) \times \pi_{r} (G_{ p_{3} } ) \times \cdots $. \qed 
\end{Lem}

\noindent
{\it Remark.}
In the recent work \cite{AS}, 
a subgroup tower of Abelian groups is introduced and studied.
A subgroup tower of a group $G$ is defined as a pair of groups 
$ (G_{1}, G_{2}, G_{3}, \cdots ) $ where $ G_{1} = G $ and
$ G_{1} \supseteq G_{2} \supseteq G_{3} \supseteq \cdots $. 
We note that definitions of a subgroup tower and a group-partition 
are similar but different. \\

Let 
\begin{align*} 
a_{r} (n) := \sum_{ |G| = n } \pi_{r} (G) ,
\end{align*}
where $G$ runs over isomorphic classes of Abelian groups of order $n$. 
We note that $ a_{1} (n) = a(n) $ since $ \pi_{1} (G) = 1 $ for any $G$. 
In order to investigate this function, 
we recall several results about the plane partition. 
Let $n$ be a positive integer and $( \lambda_{jk} )$ be a matrix 
where $ \lambda_{jk} $ are non-negative integers such that 
\begin{align*} 
\sum_{j \geq 0, k \geq 0 } \lambda_{jk} = n, \quad 
\lambda_{j1} \geq \lambda_{j2} \geq \lambda_{j3} \geq \cdots , \quad 
\lambda_{1k} \geq \lambda_{2k} \geq \lambda_{3k} \geq \cdots .
\end{align*}
Then $ (\lambda_{jk}) $ is called a plane partition of $n$. 
Let $ \mathrm{PL}_{r} (n) $ denote 
the number of plane partitions of $n$ with at most $r$-rows. 
Let $ \mathrm{PL}_{\infty} (n) $ denote 
the number of all plane partitions of $n$.
It is known  (Cf. \cite{An} p.184) that 
generating functions of the $r$-rowed plane partition $ \mathrm{PL} _{r} (n) $ 
and the plane partition $ \mathrm{PL} _{\infty} (n) $
are given by 
\begin{align} 
\sum_{ n=0 }^{ \infty } \mathrm{PL} _{r} (n) q^{n} 
= \prod_{ m=1 }^{ \infty } (1- q^{m} )^{- \min \{ m , r \} } , \qquad 
\sum_{ n=0 }^{ \infty } \mathrm{PL} _{ \infty } (n) q^{n} 
= \prod_{ m=1 }^{ \infty } (1- q^{m} )^{- m  } ,
\label{BBB}
\end{align} 
for $ |q| < 1 $. 
Now we show the following proposition which plays a key role in this paper. 

\begin{Pro} 
Let $p$ be a prime number and $n$ be a positive integer. 
For $ r = 1, 2, 3, \cdots  $, or $ r = \infty $, we have 
\begin{align} 
a_{r} (p^{n}) = \sum_{ |G| = p^{n} } \pi_{r} (G) = \mathrm{PL}_{r} (n) . 
\label{BBC} %%%%%%
\end{align} 
\end{Pro} 

\begin{proof}
We show that there is a one-to-one correspondence between 
group-partitions of Abelian $p$-groups with order $p^{n}$ 
and plane partitions of $n$. 
First we suppose that 
$ ( \lambda_{jk} ) $ 
is a $r$-rowed plane partition of $n$, that is
\begin{align*} 
\sum_{j=1}^{r} \sum_{k=1}^{n} \lambda_{jk} = n , \quad 
\lambda_{j1} \geq \lambda_{j2} \geq \cdots \geq \lambda_{jn} \geq 0 , 
\quad 
\lambda_{1k} \geq \lambda_{2k} \geq \cdots \geq \lambda_{rk} \geq 0 .
\end{align*}
We put 
\begin{align*} 
G_{j} := 
( \mathbb{Z}/ p \mathbb{Z} )^{ \lambda_{j1} - \lambda_{j2} } 
\oplus ( \mathbb{Z}/ p^{2} \mathbb{Z} )^{ \lambda_{j2} - \lambda_{j3} } 
\oplus \cdots
\oplus ( \mathbb{Z}/ p^{n} \mathbb{Z} )^{ \lambda_{jn} } ,
\end{align*} 
and $ G := G_{1} \oplus G_{2} \oplus \cdots \oplus G_{r}  $. 
Then $ G_{j} $, $G$ are Abelian $p$-groups 
such that $ | G | = p^{n} $ and $ (G_{j}) \vdash G $. 
Actually, the condition 
$ \lambda_{1k} \geq \lambda_{2k} \geq \cdots \geq \lambda_{rk} $ 
means that 
$ G_{1} \supseteq G_{2} \supseteq G_{3} 
\supseteq \cdots \supseteq G_{r} $. 

Conversely, we assume that $G$ is a Abelian group of order $ p^{n} $,  
and $ ( G_{j} ) $ is a group-partition of $G$ with at most $r$ factors. 
Let $ \mu_{jk} $ by 
\begin{align*} 
G_{j} = 
( \mathbb{Z}/ p \mathbb{Z} )^{ \mu_{j1} } 
\oplus ( \mathbb{Z}/ p^{2} \mathbb{Z} )^{ \mu_{j2} } \oplus 
\cdots
\oplus ( \mathbb{Z}/ p^{n} \mathbb{Z} )^{ \mu_{jn} } ,
\end{align*} 
and $ \lambda_{ jk } := \sum_{m=k} ^{n} \mu_{jm} $. 
Then we see that 
$ \lambda_{j1} \geq \lambda_{j2} \geq \cdots \geq \lambda_{jn} \geq 0 $. 
Since we have 
$ |G| = p^{n} $, it follows that 
$ \sum_{j=1}^{r} \sum_{k=1}^{n} \lambda_{jk} = n $.   
Moreover 
$ G_{1} \supseteq G_{2} \supseteq G_{3} 
\supseteq \cdots \supseteq G_{r} $ 
means that 
$ \lambda_{1k} \geq \lambda_{2k} \geq \cdots \geq \lambda_{rk} $. 
Hence $ (\lambda _{jk}) $ is a plane partition 
of $n$ with at most $r$ rows. 
This completes the proof of the proposition. 
\end{proof} 

\begin{Cor}
The function $ a_{r} (n) $ is multiplicative. 
\end{Cor}

\begin{proof} 
Let $ n = p_{1}^{ m_{1} } p_{2} ^{ m_{2} }  p_{3} ^{ m_{3} } \cdots  $. 
Using Lemma 2.2 and (\ref{BBC}), we have
\begin{align*} 
a_{r} (n) 
& = \sum_{ |G| 
= p_{1}^{ m_{1} } p_{2} ^{ m_{2} }  p_{3} ^{ m_{3} } \cdots } \pi_{r} (G) 
= \sum_{ |G_{1}| = p_{1} ^{ m_{1} }} \sum_{ |G_{2}| = p_{2} ^{ m_{2} }} 
\cdots \pi_{r} ( G_{1} \oplus G_{2} \oplus \cdots ) \\
& = \Big( \sum_{ |G_{1}| = p_{1} ^{ m_{1} } } \pi_{r} ( G_{1}) \Big) 
\times \Big( 
\sum_{ |G_{2}| = p_{2} ^{ m_{2} } } \pi_{r} ( G_{2}) 
\Big) \times \cdots 
= \mathrm{PL}_{r} (m_{1}) 
\times \mathrm{PL}_{r} (m_{2}) \times \cdots \\
& = a_{r} (p_{1}^{m_{1}} ) 
\cdot a_{r} (p_{2}^{m_{2}} ) 
\cdot a_{r} (p_{3}^{m_{3}} ) \cdot \cdots .
\end{align*}
This shows that $ a_{r} (n) $ is multiplicative.  
\end{proof}

\begin{Cor} 
For any $ \epsilon > 0 $, we have 
\begin{align} 
a _{r} (n) \leq 
a_{\infty} (n) 
= O _{ \epsilon } (n ^ {\epsilon } ) . 
\label{BBD} %%%%%%
\end{align} 
\end{Cor} 

\begin{proof} 
It is known (Cf. \cite{An} p.199) that 
\begin{align*}
\mathrm{PL}_{ \infty }(n) 
\sim ( \zeta(3)^{7} \cdot 2^{-11} )^{1/36} \cdot n^{-25/36} \cdot
\exp \{ 3 \cdot 2^{-2/3} \cdot \zeta(3)^{1/3} \cdot n^{2/3} +2c \} 
\end{align*}
as $ n \to \infty $. 
Here 
\begin{align*} 
c = \int_{0}^{\infty} \frac{ y \log y }{e^{ 2 \pi y } -1 } dy 
\end{align*} 
and $ \zeta (s) $ is the Riemann zeta function. 
Using this, 
we estimate roughly as 
$ \mathrm{PL}_{ \infty }(n) \leq \exp\{ C n^{2/3} \} $ 
where $C$ is a positive number. 
Then for $ \epsilon > 0 $ and 
$ n = p_{1}^{ m_{1} } p_{2} ^{ m_{2} }  p_{3} ^{ m_{3} }\cdots  $, 
we have 
\begin{align*}
0 \leq
\frac{ a _{\infty} (n) }{ n^{\epsilon} } 
& = 
\prod_{ p_{i} \leq e^{ C / \epsilon } } 
\frac{ \mathrm{PL}_{\infty} (m_{i}) }{ p_{i}^{m_{i} \epsilon } }
\times 
\prod_{ p_{i} > e^{ C / \epsilon } } 
\frac{ \mathrm{PL}_{\infty} (m_{i}) }{ p_{i}^{m_{i} \epsilon } } \\
& \leq 
\prod_{ p_{i} \leq e^{ C / \epsilon } } 
\frac{ \exp ( C m_{i}^{ 2/3 } ) }{ p_{i}^{m_{i} \epsilon } }
\times 
\prod_{ p_{i} > e^{ C / \epsilon } } 
\frac{ \exp ( C m_{i}^{ 2/3 } ) }{ \exp ( C m_{i} ) } \\
& \leq 
\prod_{ p_{i} \leq e^{ C / \epsilon } } 
\frac{ \exp ( C m_{i}^{ 2/3 } ) }{ p_{i}^{m_{i} \epsilon } } 
\leq D_{\epsilon} .
\end{align*}
Here $ D_{\epsilon} > 0 $ depends only on $ \epsilon$. 
Hence the result follows. 
\end{proof}

%%%%%%%%%%%%%%%%%%%%%%%%%%%%%%%%%%%%%%%%%%%%%%%%%%%%%%%%%%%%%%%%%%%%%%%%%%%%%%%%
%%%%%%%%%%%%%%%%%%%%%%%%%%%%%%%%%%%%%%%%%%%%%%%%%%%%%%%%%%%%%%%%%%%%%%%%%%%%%%%%

\section{Dirichlet series and Euler products } 

We first recall the classical Dirichlet $L$-function. 
Let $ \chi $ be a Dirichlet character modulo $j$ 
and $ L (s, \chi ) := \sum_{n=1}^{\infty}\chi(n) n^{-s} $ 
be the Dirichlet $L$-function. 
Then $ L( s, \chi ) $ is defined for $ \mathrm{Re} (s) > 1 $, 
and has the Euler product expression: 
$ L(s, \chi ) = \prod_{p} ( 1 - \chi (p) p^{-s} ) ^{-1} $. 
If $ \chi_{0} $ is the principal character modulo $j$, 
then the Dirichlet $L$-function $ L (s, \chi_{0} ) $ is 
essentially the Riemann zeta function: 
\begin{align*} 
L( s , \chi_{0} ) = \prod_{ p|j } ( 1 - p^{-s} ) \times \zeta (s) .
\end{align*}  
For a non-principal character $ \chi \neq \chi_{0} $, 
it is known that 
the Drichlet $L$-function $ L( s, \chi ) $ 
can be analytically continued as an entire function. 

We next consider the Dirichlet series: 
\begin{align*} 
J_{r} (s , \chi ) 
:= \sum_{ n = 1 } ^{ \infty } 
\chi (n) a_{r} (n) n^{-s} \qquad ( r = 1,2,3, \cdots , \infty ) .
\end{align*} 
We see from (\ref{BBD}) that 
this series $ J_{r} (s , \chi ) $ converges absolutely 
for $ \mathrm{Re} (s) > 1 $. 
Further, we show that 
$ J_{r} ( s, \chi ) $ has the Euler product expression as follows.

\begin{Pro} 
For $ r = 1, 2, 3, \cdots $, or $ r = \infty $, 
the function $ J_{r} (s , \chi ) $ has the Euler product expression: 
\begin{align} 
& J_{r} ( s , \chi ) 
= \prod_{ m=1 }^{ \infty } \prod_{ p : \; \mathrm{prime} } 
( 1 - \chi (p) ^{m} p^{ -m s } )^{ - \min \{ m , r \} } 
= \prod_{ m=1 }^{ \infty } L ( ms , \chi ^{m} ) ^{ \min \{ m , r \} } ,
\label{CCC} \\
& J_{\infty} ( s , \chi ) 
= \prod_{ m=1 }^{ \infty } \prod_{ p : \; \mathrm{prime} } 
( 1 - \chi (p) ^{m} p^{ -m s } )^{ -m } 
= \prod_{ m=1 }^{ \infty } L ( ms , \chi ^{m} ) ^{ m }  . \nonumber
\end{align} 
Here the product converges absolutely for $ \mathrm{Re} (s) > 1 $. 
\end{Pro} 

\begin{proof} 
Since the sum
$ \sum_{ m=1 }^{\infty } \sum_{p } m p^{-m \sigma } $ 
converges for $ \sigma > 1 $, 
it follows that 
the Euler product of (\ref{CCC})  
converges absolutely for $ \mathrm{Re} (s) > 1 $. 
For a fixed prime number $p$, 
we see from (\ref{BBB}) and (\ref{BBC}) that 
\begin{align*} 
\prod_{ m=1 }^{ \infty } 
( 1 - \chi (p) ^{m} p^{ -m s } )^{ - \min \{ m , r \} }
= \sum_{ n=0 }^{ \infty } 
\mathrm{PL}_{r} (n) \chi (p)^{n} p^{-ns} 
= \sum_{ n=0 }^{ \infty } 
a_{r} (p^{n}) \chi (p^{n}) p^{-ns} . 
\end{align*} 
Since $a_{r} (n)$ is a multiplicative function, 
we obtain (\ref{CCC}).  
\end{proof} 

\begin{Cor}
The function $ J_{r} ( s, \chi ) $ defined for $ \mathrm{Re} (s) > 1 $ 
has a meromorphic continuation to the half plane 
$ \mathrm{Re} (s) > 0 $.  
\end{Cor}

\begin{proof}
Using Proposition 3.1, we have 
\begin{align} 
J_{r} (s , \chi )
= \prod_{ m=1 }^{ M-1 } L ( ms , \chi ^{m} ) ^{ \min \{ m , r \} } 
\times \prod_{ m=M }^{ \infty } \prod_{ p } 
( 1 - \chi (p) ^{m} p^{ -m s } )^{ - \min \{ m , r \} } , 
\label{CCF}
\end{align} 
for $ \mathrm{Re} (s) > 1 $. 
Since the sum
$ \sum_{ m=M }^{\infty } \sum_{p } m p^{-m \sigma } $ 
converges for $ \sigma > 1/M $, 
it follows that 
the second term of (\ref{CCF})  
can be defined for $ \mathrm{Re} (s) > 1/M $ for any $M$. 
This gives a meromorphic continuation of $ J_{r} (s, \chi) $
to $ \mathrm{Re} (s) > 0 $. 
\end{proof} 

\noindent
{\it Remark.}
Using the result in \cite{Da}, 
we see that 
the imaginary axis is a natural boundary 
of the function $ J_{r} (s, \chi) $. \\

%%%%%%%%%%%%%%%%%%%%%%%%%%%%%%%%%%%%%%%%%%%%%%%%%%%%%%%%%%%%%%%%%%%%%%%%%%%%%%%%
%%%%%%%%%%%%%%%%%%%%%%%%%%%%%%%%%%%%%%%%%%%%%%%%%%%%%%%%%%%%%%%%%%%%%%%%%%%%%%%%

\section{Asymptotic estimates of group-partitions}

To describe the distribution of group-partitions, 
we recall two lemmas. 

\begin{Lem}[Cf. \cite{Mu} p.26] 
For positive integers $j$ and $a$ with $ (j,a) = 1 $, we have
\begin{align*} 
\frac{1}{ \phi (j) } \sum_{ \chi \bmod j } 
\bar{ \chi } (a) \chi (b) 
= \begin{cases}
1 & \textit{if} \quad b \equiv a \; ( \bmod \, j ) , \\
0 & \textit{if} \quad b \equiv a \; ( \bmod \, j ) , 
\end{cases}
\end{align*} 
where $ \chi $ runs 
over the distinct Dirichlet characters modulo $j$, 
and 
$ \phi (j) = j \prod_{ p | j } ( 1-p^{-1} ) $ 
is the Euler function. \qed 
\end{Lem}  

The following lemma is known as the Ikehara-Wiener theorem. 

\begin{Lem}[ Cf. \cite{Mu} p.43]
Let $ F(s) : = \sum_{n=1}^{\infty} b_{n} n^{-s} $ 
be the Dirichlet sereis 
with positive real coefficients 
and absolutely convergent for $ \mathrm{Re} (s) > 1 $. 
Suppose that $ F(s) $ can be extended to a meromrphic function 
in the region $ \mathrm{Re} (s) \geq 1 $ having no poles 
except for a simple pole at $ s = 1 $ 
with residue $ R \geq 0 $. 
Then 
\begin{align*} 
\sum_{ n \leq x } b_{n} = R x + o(x)
\end{align*} 
as $ x \to \infty $. \qed
\end{Lem}  

We now prove the main theorem. 

\begin{The} 
For $ r = 1, 2, 3, \cdots $, or $ r = \infty $, we have 
\begin{align*} 
\sum_{ n \leq x } a_{r} (n)
= \Big\{ \prod_{ m=2 }^{ \infty } 
\zeta (m) ^{ \min \{ m , r \} } \Big\} \; x 
+ o(x) 
\qquad \text{as } x \to \infty . 
\end{align*} 
More precisely, 
if $j$ and $k$ are positive integers with $ (j,k) = 1 $, 
then we have 
\begin{align} 
\sum_{ 
\substack{ n \leq x  \\ 
n \equiv k \;( \bmod \, j ) } } 
a _{r} (n) 
= \frac{ c_{r} (j) }{ \phi (j) } x 
+ o(x) \qquad \text{as } x \to \infty ,
\label{FFF}
\end{align} 
where $ \phi (j) $ denotes the Euler function 
and $ c _{r} (j) $ is a constant which is given by 
\begin{align} 
c _{r} (j)  
:= \prod_{ m=1 }^{ \infty } \prod_{ p|j } 
( 1-p^{-m} )^{ \min \{ m , r \} } 
\times \prod_{ m=2 }^{ \infty } \zeta (m) ^{ \min \{ m , r \} } .
\label{EEF}
\end{align} 
\end{The}

\begin{proof}
For $ \mathrm{Re} (s) > 1 $, we consider the characteristic function: 
\begin{align*} 
C _{r} (s,k) 
:= \sum_{ 
\substack{ n =0  \\ 
n \equiv k \;( \bmod \, j ) } } ^{\infty}
a _{r} (n) n^{-s} . 
\end{align*} 
From Lemma 4.1, we have  
\begin{align*} 
C _{r} (s,k) 
=& \sum_{ n=0 }^{ \infty } 
\Big\{ 
\frac{1}{ \phi (j) } 
\sum_{ \chi } \bar{\chi} (k) \chi (n) 
\Big\} a_{r} (n) n^{-s}
= \frac{1}{ \phi (j) } \sum_{ \chi }
\bar{ \chi } (k) J_{r} (s , \chi ) \\
=& \frac{1}{ \phi (j) } J _{r} ( s , \chi_{0} )
+ \frac{1}{ \phi (j) } \sum_{ \chi \neq \chi_{0} }
\bar{ \chi } (k) J_{r} (s , \chi ) . 
\end{align*} 
We see that 
$ C _{r} (s,k) $ has an analytic continuation 
as a meromorphic function to $ \mathrm{Re} (s) > 0 $ 
with a simple pole at $ s=1 $. 
The residue of $ C _{r} (s,k) $ at $ s = 1 $ is given by 
\begin{align*} 
\lim_{ s \to 1 } (s-1) C _{r} (s,k) 
&= \frac{1}{ \phi(j) } \lim_{s \to 1 } (s-1) J_{r} ( s , \chi_{0} ) \\
&= \frac{1}{ \phi(j) } 
\prod_{ m=1 }^{ \infty } \prod_{ p|j } ( 1-p^{-m} ) ^{ \min \{ m , r \} } 
\times \prod_{ m=2 }^{ \infty } \zeta (m) ^{ \min \{ m , r \} } 
= \frac{ c_{r} (j) }{ \phi(j) }. 
\end{align*} 
Applying Lemma 4.2 
to the function $ C _{r} (s,k) $, 
we obtain the theorem. 
\end{proof} 

For an isomorphic class $G $ of finite Abelian groups, 
we put $ \sigma _{r} (G) $ by the number of group-partitions of $G$ 
with just $r$ factors. 
Then we see that 
$ \sigma_{r} (G) = \pi_{r} (G) - \pi_{r-1} (G)$. 
Let 
$ b_{r} (n) := \sum_{|G|=n} \sigma_{r} (G) $. 
Then, for $ r = 2, 3, 4, \cdots $, it is clear that 
\begin{align*} 
\sum_{ 
\substack{ n \leq x  \\ 
n \equiv k \;( \bmod \, j ) } } 
b _{r} (n) 
= \frac{ c_{r} (j) - c_{r-1} (j) }{ \phi (j) } x  
+ o(x) \qquad \text{as } x \to \infty 
\end{align*} 
where $ c _{r} (j) $ is a constant given by (\ref{EEF}) .\\

\noindent
{\it Remark.}
The error term of (\ref{AAB}) are studied by several authors 
(Cf. \cite{ES},\cite{Li} for $j=1$, 
\cite{Iv},\cite{Sa} for $ j \in \mathbb{Z}_{ \geq 1 } $ ). 
Estimations of error terms in (\ref{FFF}) 
will be studied in somewhere else. \\

%%%%%%%%%%%%%%%%%%%%%%%%%%%%%%%%%%%%%%%%%%%%%%%%%%%%%%%%%%%%%%%%%%%%%%%%%%%%%%%%
%%%%%%%%%%%%%%%%%%%%%%%%%%%%%%%%%%%%%%%%%%%%%%%%%%%%%%%%%%%%%%%%%%%%%%%%%%%%%%%%

Graduate School of Mathematics, Kyushu University 

6-10-1, Hakozaki Fukuoka 

812-8581, Japan

E-mail momo@math.kyushu-u.ac.jp

\end{document}